\documentclass[10pt,oneside,reqno]{amsart}

\usepackage{eucal, amssymb,amsmath,graphicx, amsfonts, latexsym}
\usepackage[a4paper,left=1.5cm,top=2cm,bottom=2cm,right=2cm,
  heightrounded]{geometry}
\usepackage[utf8]{inputenc}
\usepackage[all,arc]{xy}
\usepackage{enumerate}
\usepackage{mathrsfs}
\usepackage{verbatim}
\usepackage{lmodern}
\usepackage{float}
\usepackage{setspace}
\usepackage{textcomp}
\usepackage{subcaption}
\usepackage{mathtools}
\usepackage[outdir=./]{epstopdf}
\usepackage{setspace} 
\theoremstyle{definition}

\theoremstyle{remark}

\makeatletter
\let\c@equation\c@thm
\makeatother
\numberwithin{equation}{section}

\title{O\MakeLowercase{n the symplectic integration of the} K\MakeLowercase{lein} 
G\MakeLowercase{ordon} \MakeLowercase{lattice model}}

\date{}
\author{B. S\MakeLowercase{enyange}\\
D\MakeLowercase{epartment of} M\MakeLowercase{athematics and} 
A\MakeLowercase{pplied} M\MakeLowercase{athematics}\\
U\MakeLowercase{niversity of} C\MakeLowercase{ape} T\MakeLowercase{own}, 
R\MakeLowercase{ondebosch, 7701}, S\MakeLowercase{outh} A\MakeLowercase{frica}}

\begin{document}
\onehalfspacing
\begin{abstract}

We investigate the performance of various methods of symplectic 
integration, which are based on two part splitting of the integration 
operator, for the numerical integration of multidimensional Hamiltonian 
systems. We implement these schemes to study the behaviour of the 
one-dimensional quartic Klein Gordon disordered lattice with many 
degrees of freedom (of the order of a few hundreds) and compare their 
efficiency for the weak chaos regime of the system’s dynamics. For this 
reason we perform extensive numerical simulations for each considered 
integration scheme. In this process, the second moment and the 
participation number of propagating wave packets, along with the system’s 
relative energy error and the required CPU time are registered and compared.

\vspace{0.1in}

\end{abstract}

\linespread{0.5}
\maketitle
\begin{enumerate}
 \item[{\bf 1.}] {\bf Introduction}\\
Symplectic integrators (SIs) are known to preserve the symplectic nature of the Hamiltonian system 
and keep bounded the error of the computed value of the Hamiltonian. This is one of the advantages that these integrators have over general purpose 
integrators. In (\cite{2},\cite{22},\cite{4}) various SIs have been applied in the study of the chaotic behavior of 
two one-dimension Hamiltonian lattices, namely the Klein-Gordon (KG) chain and the Discrete NonLinear 
Schr$\ddot{o}$dinger model. These studies showed that there exist different dynamical behaviors, namely 
the so called weak chaos, strong chaos and the self trapping regime. In this study we consider a wider range of 
SIs for integration of multidimensional Hamiltonian systems.\\ 
In the next section we give a brief discussion of the KG lattice as the Hamiltonian model 
to use in this study. We also give a description of SIs of generalised order, order 
two and order four with an insight of how composition techniques are used to generate schemes of higher 
order. Section $3$ is devoted to comparing the performance of these SIs for the integration of the 
KG lattice after which we present our conclusions in section $4$.

 \item[{\bf 2.}] {\bf The KG Hamiltonian model and the integration schemes}
 
The Hamiltonian $H$ of the one-dimensional KG lattice model of coupled anharmonic oscillators with $N$ degrees of freeedom is 
\begin{equation}\label{ham}
 H({\bf q},{\bf p}) = \sum_i \frac{p_i^2}{2} + \frac{\epsilon_i}{2}q_i^2 + \frac{1}{4}q_i^4 + \frac{1}{2W}\left(q_{i+1} - q_i\right)^2,
\end{equation}

where $q_i$ and $p_i$ are the generalised position and momenta of site $i$ respectively. $\epsilon_i$ are 
potential strengths that are chosen uniformly from the interval $[\frac{1}{2},\frac{3}{2}]$, and $W$ is a parameter that determines the 
extent of disorder in the lattice.
From (\ref{ham}), the resulting equations of motion

\begin{equation}\label{kgh}
  \frac{dq_i}{dt} = p_i, \qquad \qquad
 \frac{dp_i}{dt} = -\epsilon_iq_i - q_i^3 + \frac{1}{W}\left(q_{i+1} + q_{i-1} -2q_i\right) 
\end{equation} can be written as ${\frac{d{\bf z}}{dt}= \{{\bf z},H\} =: L_H{\bf z}}$, where ${{\bf z} = ({\bf q},{\bf p})}$, $L_H$ 
is the so called Poisson bracket 
that is defined by ${\{F,G\} := \sum_i\left( \frac{\partial F}{\partial q_i}\frac{\partial G}{\partial p_i} - 
\frac{\partial F}{\partial p_i}\frac{\partial G}{\partial q_i}\right)}$, for differentiable functions $F({\bf z})$ and $G({\bf z})$. 
Using initial conditions ${\bf z_0}={\bf z}(0)$, we therefore get a formal solution 
\begin{equation*}
 {\bf z}(t) = \sum_{i \geq 0}\frac{t^i}{i!}L_H^i{\bf z_0} = e^{tL_H}{\bf z_0}.
\end{equation*}
\newpage
The Hamiltonian (\ref{ham}) can be split into two integrable parts as $H({\bf z}) = A({\bf p}) + B({\bf q})$ where  
\begin{equation}
A = \sum_i \frac{p_i^2}{2},\qquad \qquad
 B = \sum_i \frac{\epsilon_i}{2}q_i^2 + \frac{1}{4}q_i^4 + \frac{1}{2W}\left(q_{i+1} - q_i\right)^2
\end{equation} and the action of the operators $e^{\tau L_A}$ and $e^{\tau L_B}$ is known analytically.


A SI approximates the operator $e^{\tau L_H}$ by a product of operators $e^{a_i\tau L_A}$ and $e^{b_i\tau L_B}$ 
where the constants $a_i$ and $b_i$ are chosen depending on the required order of the SI \cite{15}.

In our study we consider the performance of order two, order four and  generalised order SIs $ABA82,$ $ABA864,$ $ABAH864$ of \cite{23} in integrating system 
(\ref{ham}).\\
The order two SIs $SABA_2$ and $SBAB_2$ \cite{3,5} are 
\begin{equation*}
SABA_2=e^{c_1\tau L_A}e^{d_1\tau L_B}e^{c_2\tau L_A}e^{d_1\tau L_B}e^{c_1\tau L_A}
\end{equation*}
where $c_1=\frac{1}{2} - \frac{1}{2\sqrt{3}}$, $c_2 = \frac{1}{\sqrt{3}}$, $d_1 = \frac{1}{2}$, and

\begin{equation*}
SBAB_2 = e^{c_1\tau L_B}e^{d_1\tau L_A}e^{c_2\tau L_B}e^{d_1\tau L_A}e^{c_1\tau L_B}
\end{equation*} where $d_1 = \frac{1}{2}$, $c_1 = \frac{1}{6}$ and $c_2 = \frac{2}{3}$

The order of the SIs $SABA_2$ and $SBAB_2$ can be improved to order four by including a corrector term
\begin{equation*}
 {\bf C} = e^{-\tau^3\varepsilon^2\frac{c}{2}L_{\{\{A,B\},B\}}}
\end{equation*}

where the value of $c$ is $\frac{(2-\sqrt{3})}{24}$ for $SABA_2$ 
and $\frac{1}{72}$ for $SBAB_2$. We therefore get the so called $SABA_2$ with corrector, $SABA_2wc = {\bf C}SABA_2{\bf C}$ 
and $SBAB_2$ with corrector $SBAB_2wc = {\bf C}SBAB_2{\bf C}$
\newline
We also consider SIs of generalised order  $(8,2)$ and $(8,6,4)$ \cite{6,12}.

Order $(8,2)$ SI $ABA82$:
\begin{equation*}
ABA82 = e^{a_1\tau L_A}e^{b_1\tau L_B}e^{a_2\tau L_A}e^{b_2\tau L_B}e^{a_3\tau L_A}e^{b_2\tau L_B}
e^{a_2\tau L_A}e^{b_1\tau L_B}e^{a_1\tau L_A}
\end{equation*} with the coefficients $a_i$, $b_i$ stated in \cite{23} and order $(8,6,4)$ SIs
$ABA864$ and $ABAH864$ as defined in \cite{12} together with the corresponding coefficients.\newline
For a symmetric order two integrator $S_2$, an order four integrator 
$S_2Y4$ was constructed in \cite{7} by Yoshida using composition techniques. That is to say,

\begin{equation*}
 S_2Y4(\tau) := S_2(a_1\tau)S_2(a_0\tau)S_2(a_1\tau) 
\end{equation*} where

$a_0 = \frac{2^{\frac{1}{3}}}{2-2^{\frac{1}{3}}}$ and $a_1 = \frac{1}{2-2^{\frac{1}{3}}}$. 
In particular we study the behavior of order four SIs $SABA_2Y4$ 
and $SBAB_2Y4$ 

In \cite{8}, using the order two SI Leap-Frog ($LF$) an order four scheme $Sz4$ was constructed:

\begin{equation*}
 Sz4(\tau) = LF(k{\bf z})LF((1-2k){\bf z})LF(p{\bf z})
\end{equation*}
where 	. It can be easily seen that when the Yoshida composition technique is 
applied to $LF$, one ends up with $Sz4$.\newline
Forest and Ruth constructed a fourth order SI which we shall call $FRo4$ defined as 
$$FRo4 = e^{a_1\tau L_A}e^{b_1\tau L_B}e^{a_2\tau L_A}e^{b_2\tau L_B}e^{a_3\tau L_A}e^{b_3\tau L_B}
e^{a_4\tau L_A}e^{b_3\tau L_B}e^{a_3\tau L_A}e^{b_2\tau L_B}
e^{a_2\tau L_A}e^{b_1\tau L_B}e^{a_1\tau L_A}$$ with the coefficients $a_i$, $b_i$ as specified in \cite{9}.
\newpage

 \item[{\bf 3.}] {\bf Numerical results}
 
 We consider a disorder realization of $H$ in (\ref{ham}) for a total of $1000$ sites with a random 
 value of $\epsilon_i$ at site $i$. Fixing $W=4$ and $0$ initial displacement we make an initial 
 excitation of the central site with a total energy of $0.4$. We then keep track of the second moment 
 $m2$, participation number $P$ and CPU time.

 The energy of site $i$ at a time $t$ is 
 \begin{equation}    
h_i = \frac{p_i^2}{2} + \frac{\epsilon_i}{2}q_i^2 + 
 \frac{1}{4}q_i^4 + \frac{1}{4W}\left(q_{i+1} - q_i\right)^2
 \end{equation} and $\bar{i} := \sum_ii\frac{h_i}{E_t}$.
 
 With energy $E_t$ at time $t$, a normalised energy distribution $\frac{h_i}{E_t}$ of 
 site $i$, $m2 = \sum_i(i-\bar{i})^2\frac{h_i}{E_t}$ is a measure of the rate at which the wave packet spreads 
 from the initially excited central site to all sites in the lattice and $P = \frac{1}{\sum_i\left(\frac{h_i}{E_t}\right)^2}$ 
 quantifies the proportion of excited sites in the entire lattice. 
 
 In order to compare the performance of the different SIs, we adjust the time step $\tau$ so that 
 the absolute relative energy error $REe := \left|\frac{E_t - E_0}{E_0}\right| \lesssim 10^{-5}$ at a time $t$ of the evolution; where 
 $E_0$ and $E_t$ are the energies of the system at times $0$ and $t$ respectively. For each of the SIs we ensure that there 
 is a global consistance amongst the SIs in the evolution of $m2$ and $P$ for capturing the dynamics of the 
 wave packet. For each of the SIs we then record the CPU time which is required to perform the simulations.
 
 FIGURE 1 shows the results obtained when we integrate (\ref{kgh}) using order two schemes $SBAB_2$ (red curve) 
 and $SABA_2$ (green curve) and generalised order scheme $ABA82$ (gray curves). From this figure we see that 
 with all the schemes portraying practically the same dynamical behavior of the wave packet with respect to $m2$ and $P$,
 $ABA82$ has the best performance since it requires the least CPU time compared to the other SIs.
 
 In FIGURE 2 we have results for the integration of equations of motion \ref{ham} using order four schemes $SABA2wc$ 
(red curve), $Sz4$ (blue curve), $SBAB_2Y4$ (pink curve) and 
generalized order schemes $ABAH864$ (green curve) and $ABA82$ (gray curves). The generalised order scheme $ABA82$ requires 
the highest CPU time whereas the schemes $SABA2wc$ and $ABAH864$ show a better performance compared to other schemes since they require 
the least CPU time.

 \begin{figure}[htp]
    \centering
    \begin{subfigure}{0.3\textwidth}
        \includegraphics[width=\textwidth]{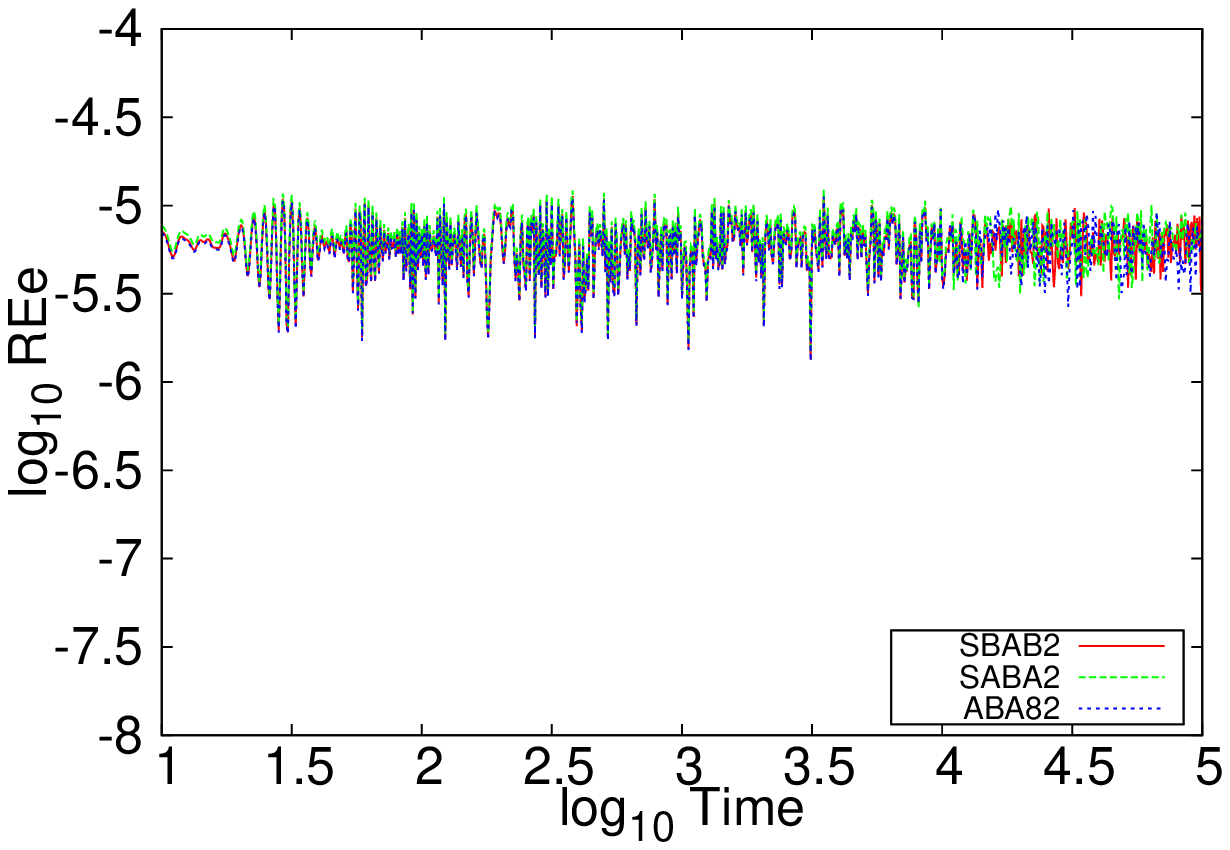}
        
    \end{subfigure}
    \begin{subfigure}{0.3\textwidth}
        \includegraphics[width=\textwidth]{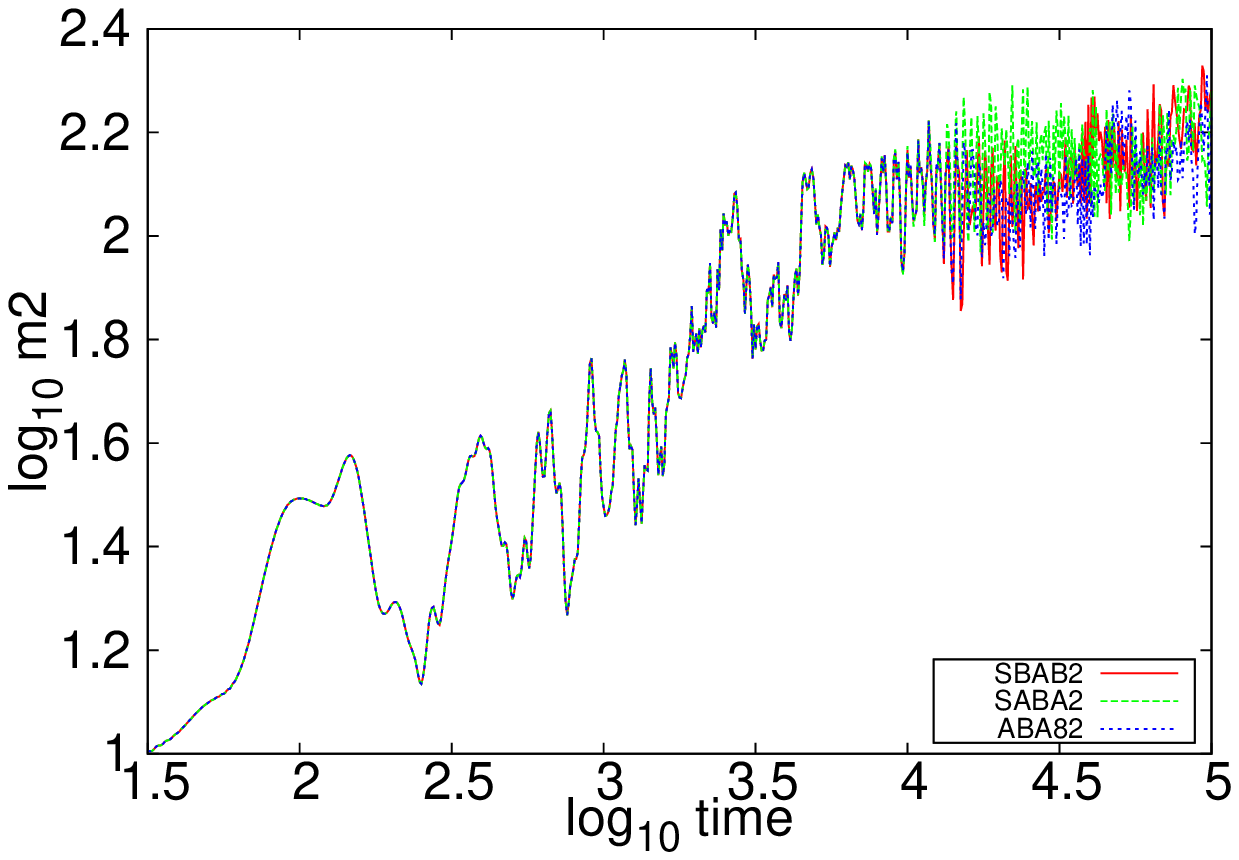}
        
    \end{subfigure}
    
    \begin{subfigure}{0.3\textwidth}
        \includegraphics[width=\textwidth]{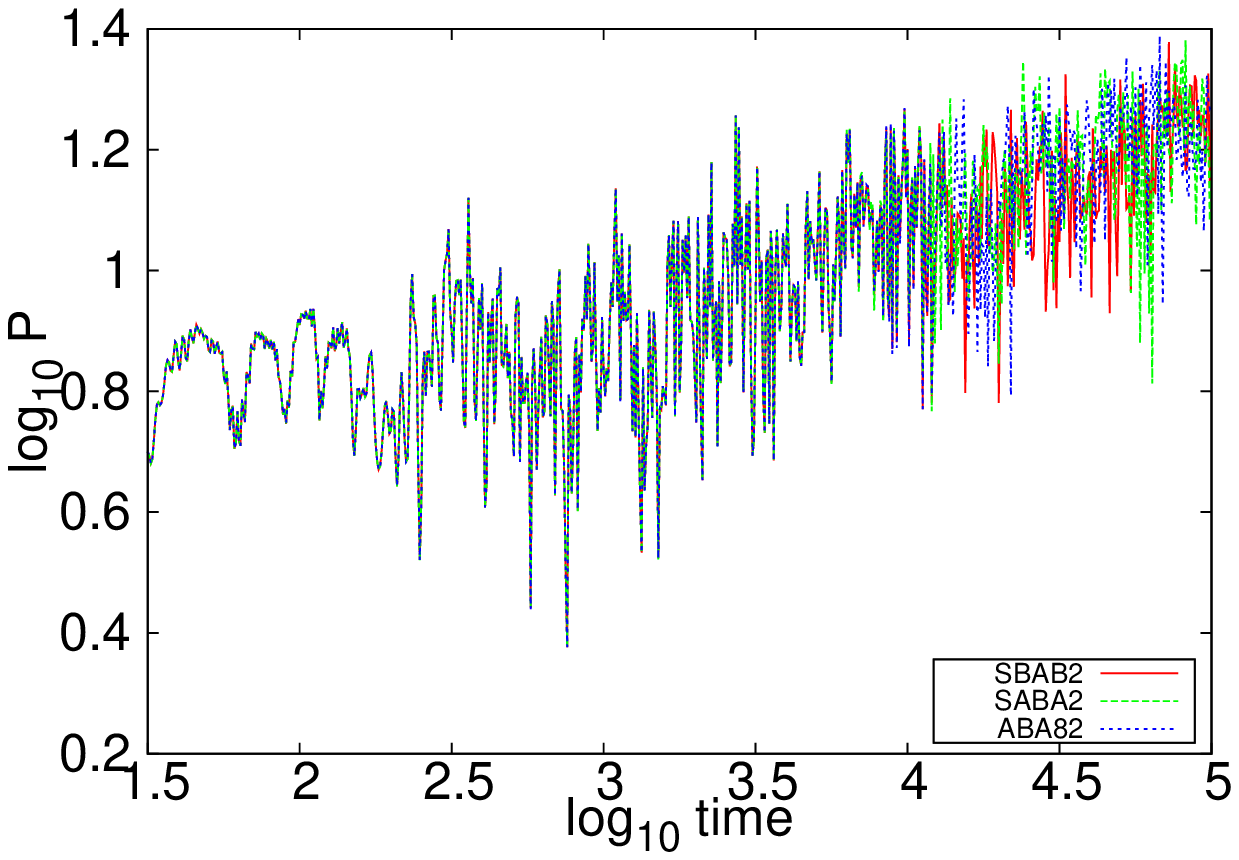}
        
    \end{subfigure}
     \begin{subfigure}{0.3\textwidth}
        \includegraphics[width=\textwidth]{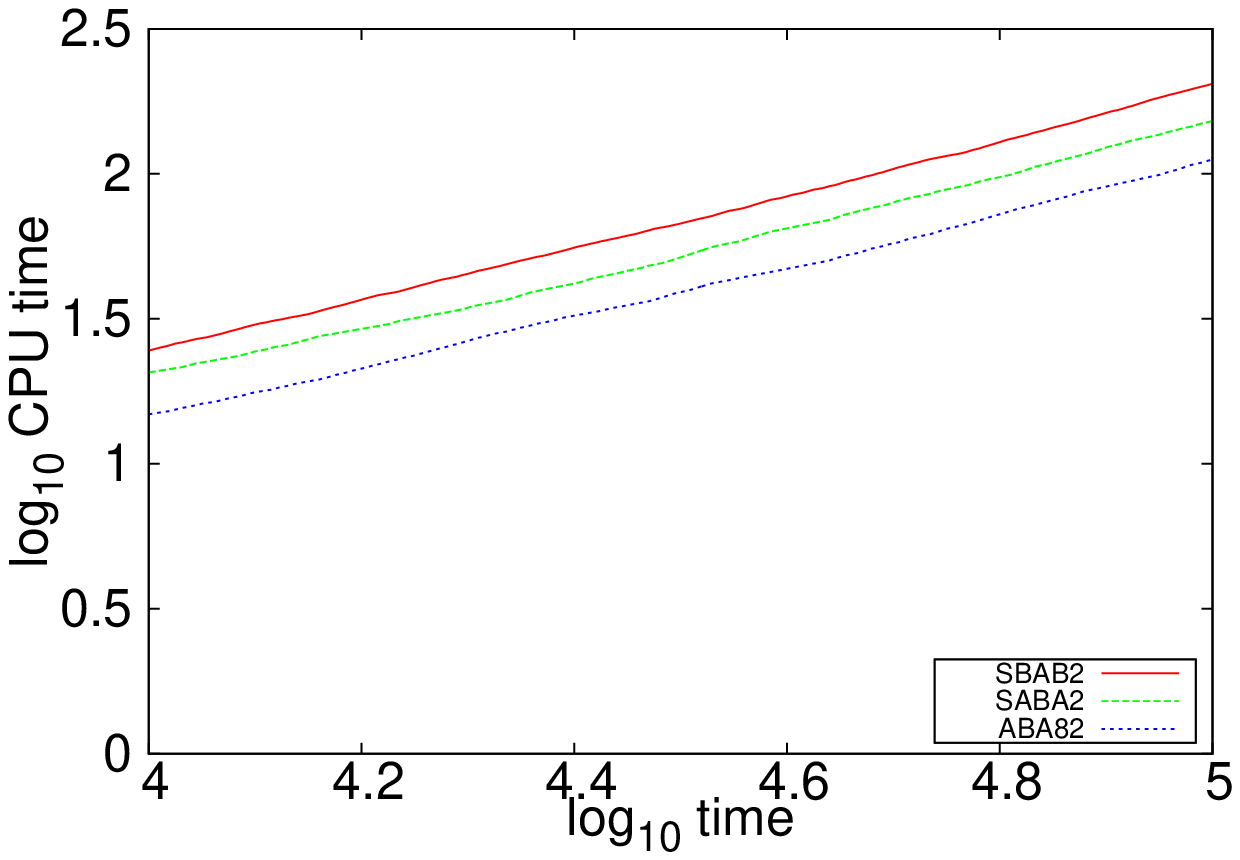}
        
    \end{subfigure}

\caption{{\tiny {\bf Results for the integration of (\ref{kgh}), by the order two schemes $SBAB_2$ for 
$\tau = 0.016$ (red curve) and $SABA_2$ for $\tau = 0.0185$ (green curve) and $ABA82$ of generalised order $(8,2)$ for $\tau = 0.032$ 
(gray curves). The panels show the logarithms of the relative energy error, second moment, participation number 
and the CPU time required for evolution upto time $10^5$.}}}
\end{figure}
\newpage
 \begin{figure}[htp]
    \centering
    \begin{subfigure}{0.3\textwidth}
        \includegraphics[width=\textwidth]{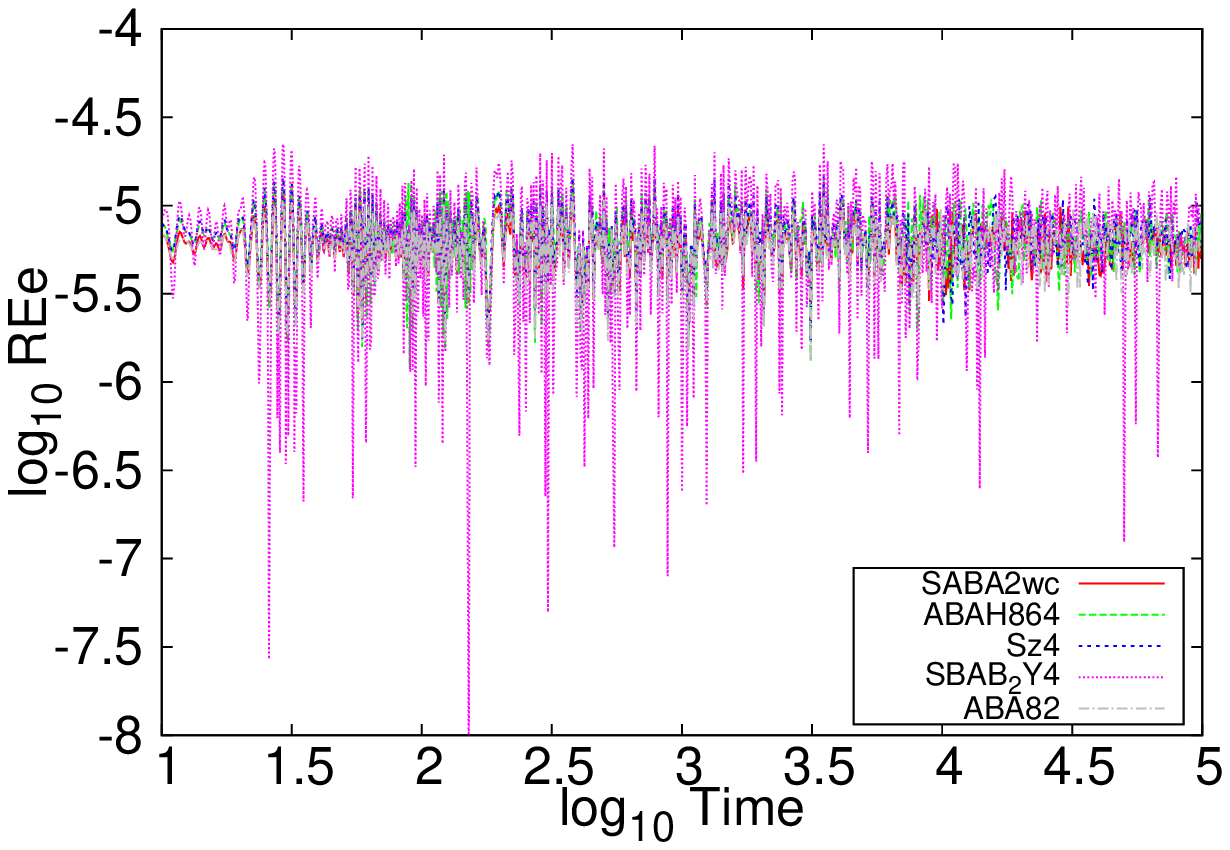}
        
    \end{subfigure}
    \begin{subfigure}{0.3\textwidth}
        \includegraphics[width=\textwidth]{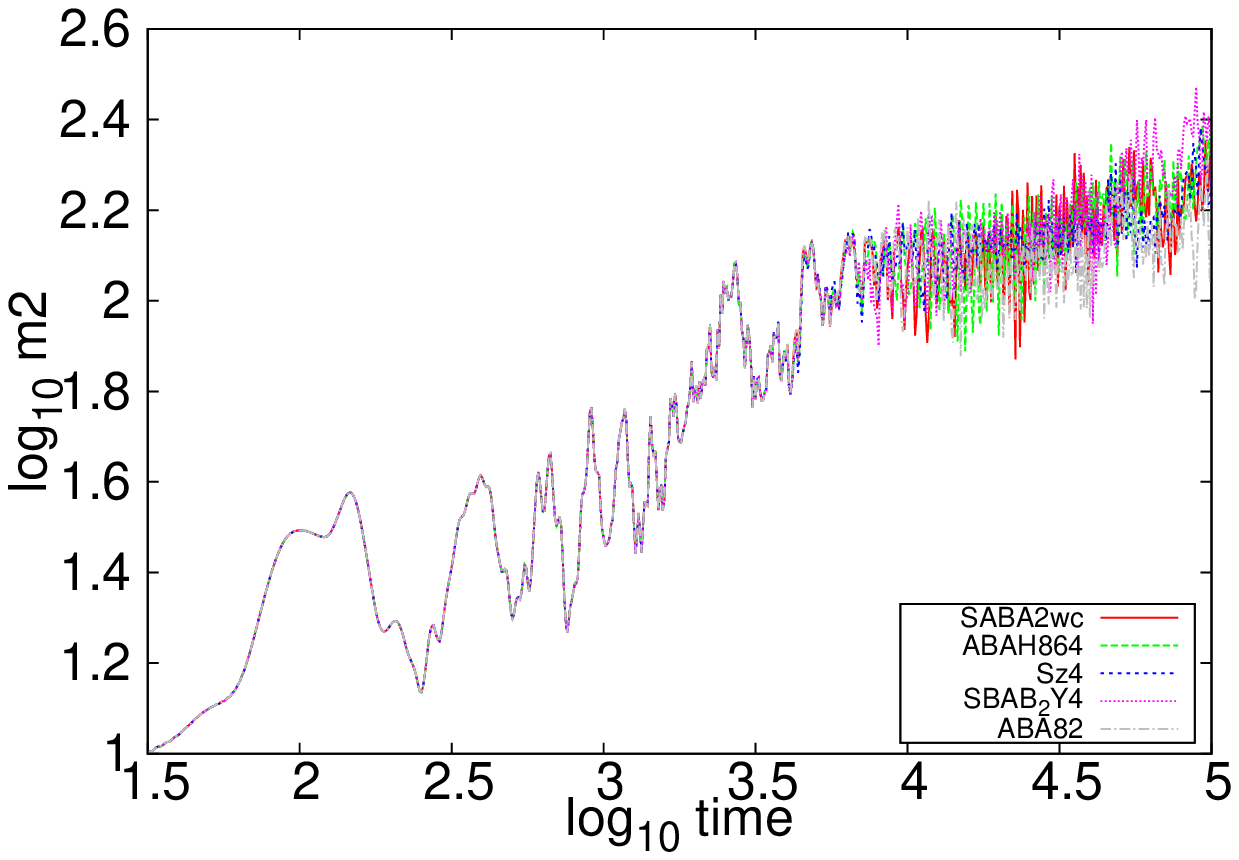}
        
    \end{subfigure}
    
    \begin{subfigure}{0.3\textwidth}
        \includegraphics[width=\textwidth]{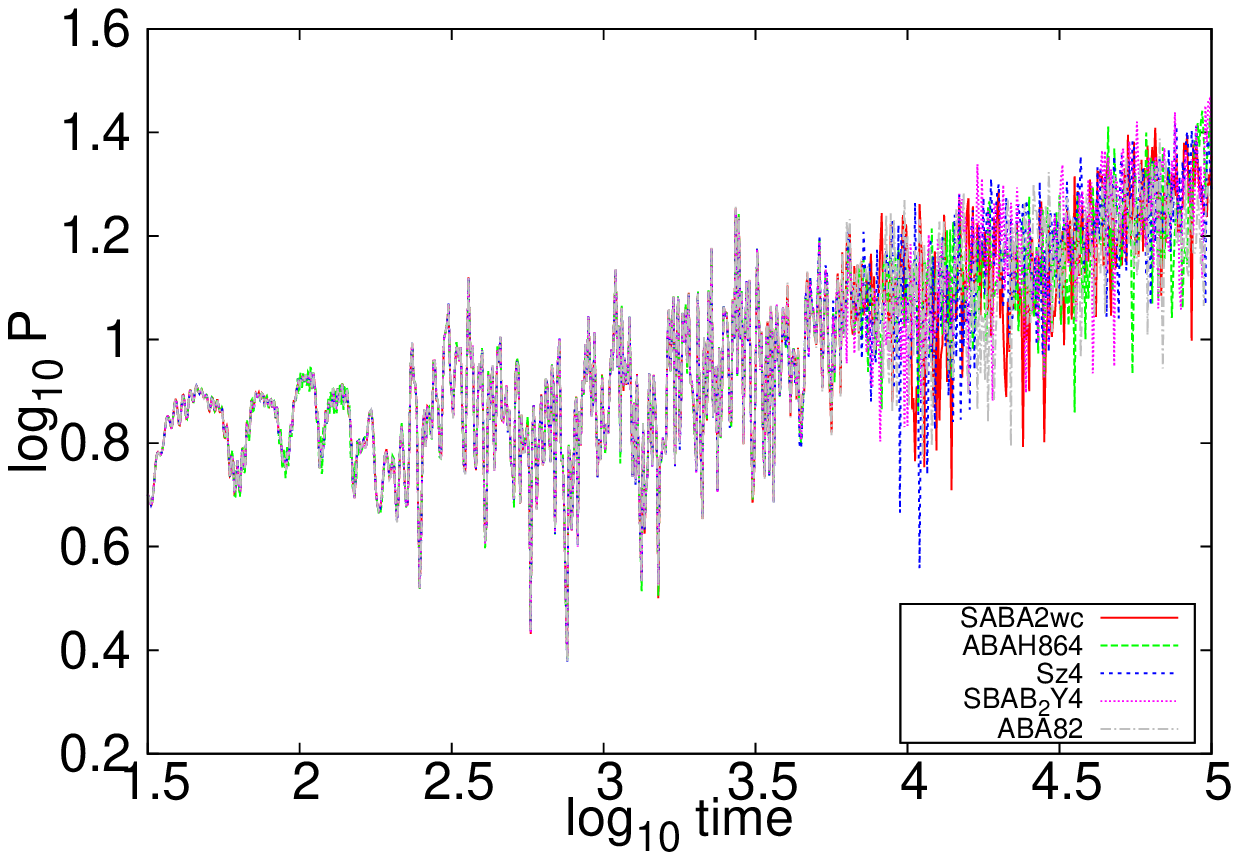}
        
    \end{subfigure}
     \begin{subfigure}{0.3\textwidth}
        \includegraphics[width=\textwidth]{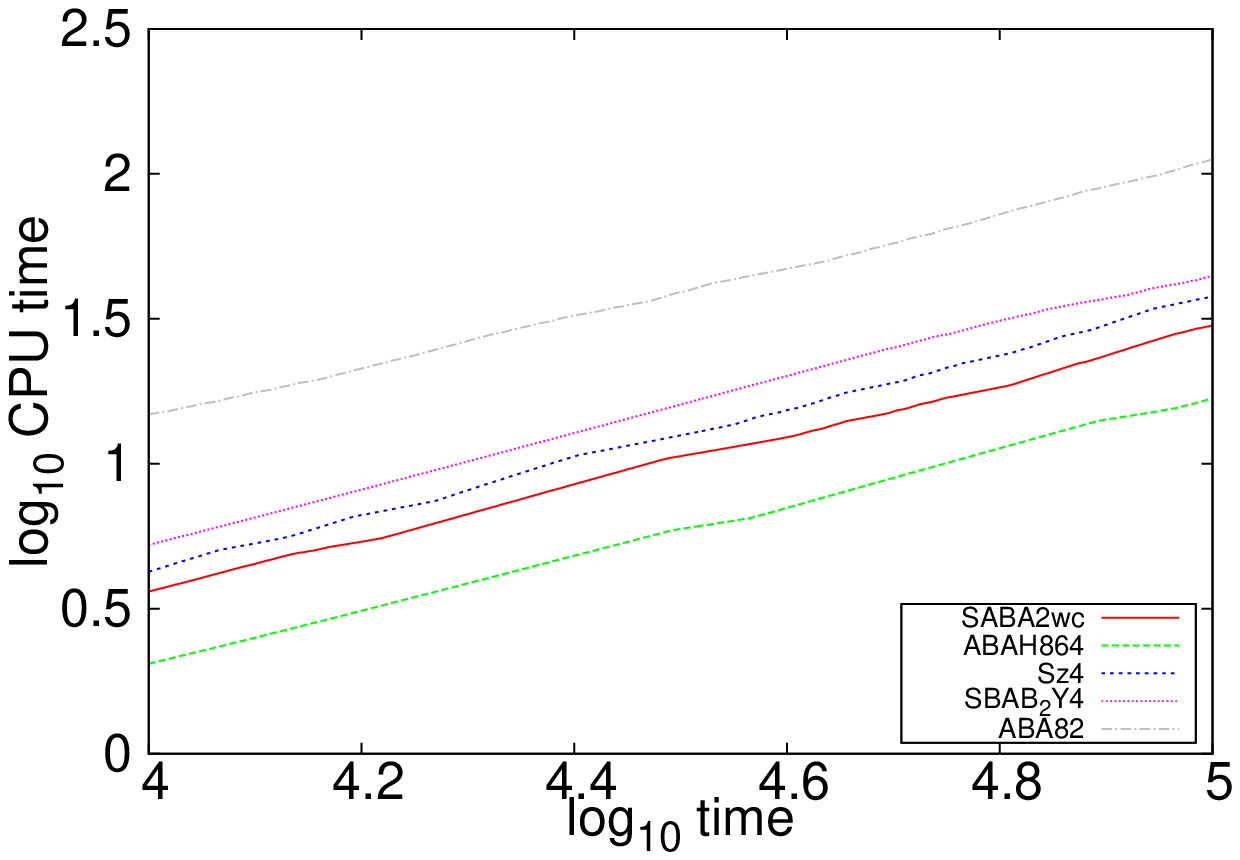}
        
    \end{subfigure}

\caption{{\tiny {\bf Results for the integration of (\ref{kgh}), by the order four schemes $SABA2wc$ for 
$\tau = 0.165$ (red curve), $ABAH864$ for $\tau = 0.355$ (green curve), $Sz4$ for $\tau = 0.084$ 
(blue curve), $SBAB_2Y4$ for $\tau = 0.13$ (pink curve) and $ABA82$ for $\tau = 0.032$ 
(gray curves).  The panels show the logarithms of the relative energy error, second moment, participation number 
and the CPU time required for evolution upto time $10^5$.}}}
\end{figure}

 In FIGURE 3 we have the results for the integration when we use $SABA_2Y4$ (red curves), 
 $SBAB_2wc$ (green curve), $ABAH864$ (gray curve), $ABA864$ (pink curve) and $FRo4$ (light blue curves).
 The SIs of generalised order $ABAH864$ and $ABA864$ show a better performance compared to all 
 the other SIs that have been used in the simulations.
 The SI $ABA864$ reveals the best performance in terms of least CPU time amongst all the SIs that were used in this work.
 
 \begin{figure}[htp]
    \centering
    \begin{subfigure}{0.3\textwidth}
        \includegraphics[width=\textwidth]{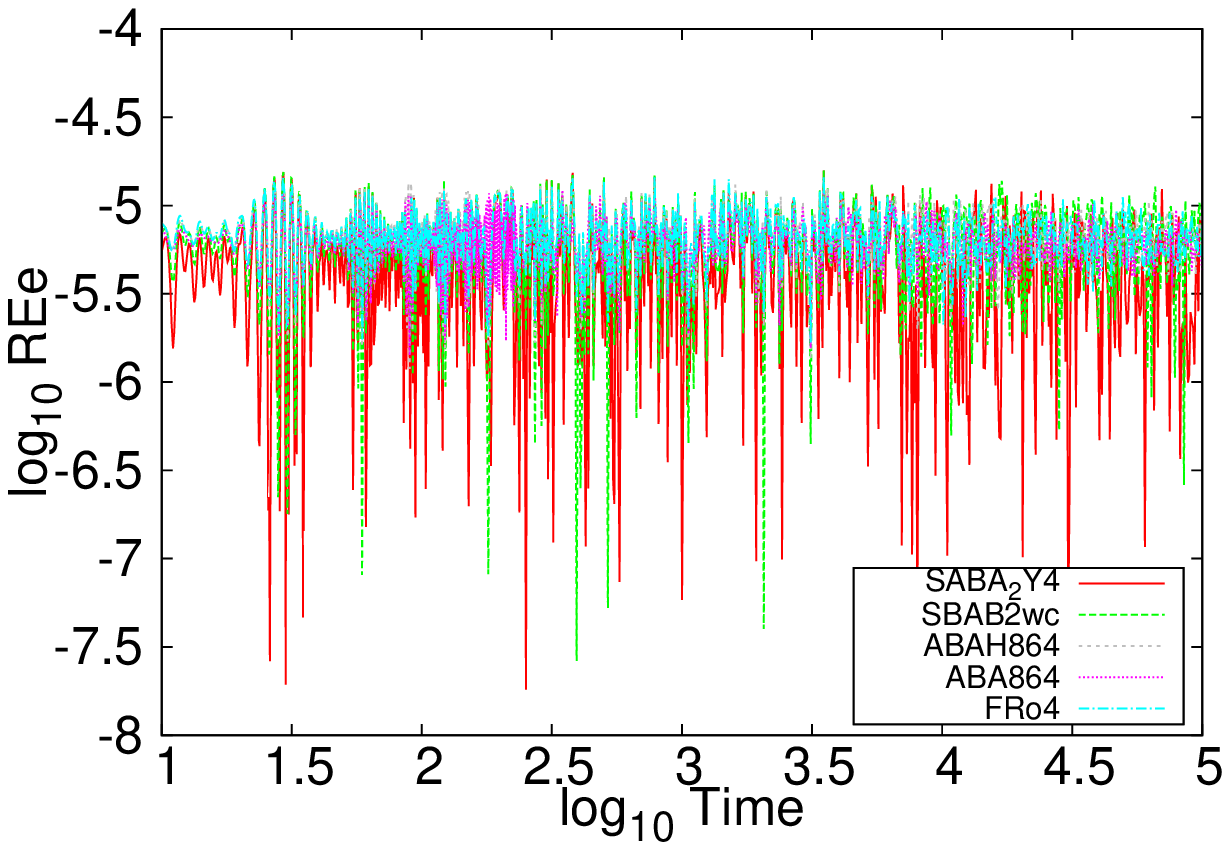}
        
    \end{subfigure}
    \begin{subfigure}{0.3\textwidth}
        \includegraphics[width=\textwidth]{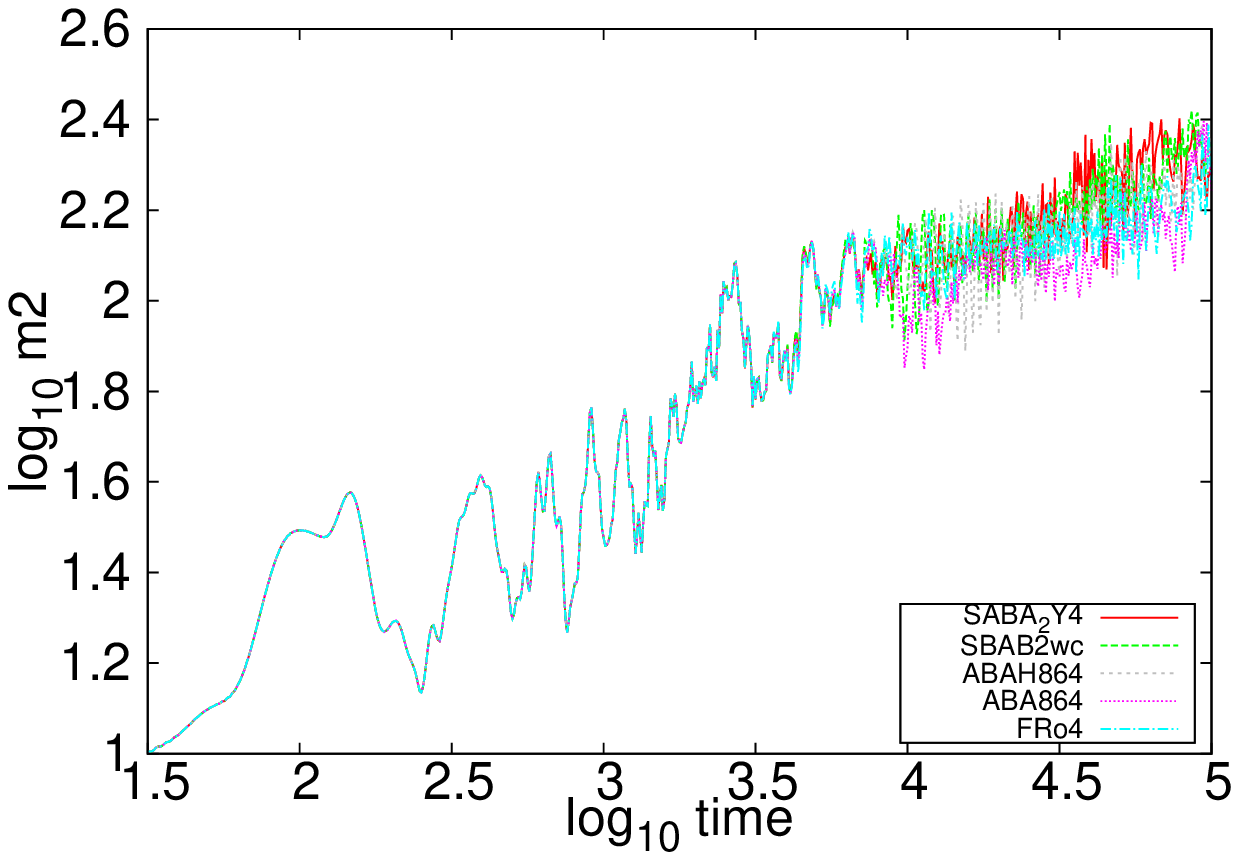}
        
    \end{subfigure}
    
    \begin{subfigure}{0.3\textwidth}
        \includegraphics[width=\textwidth]{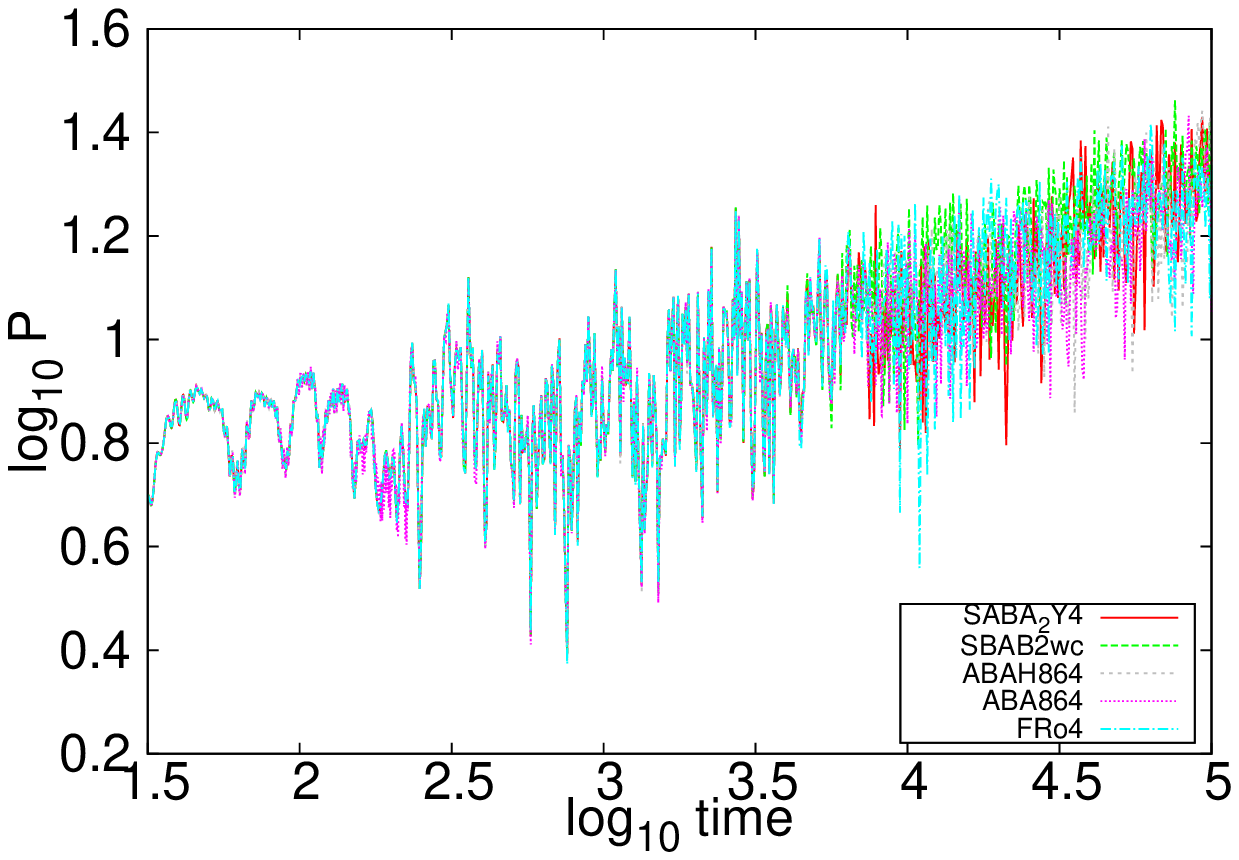}
        
    \end{subfigure}
     \begin{subfigure}{0.3\textwidth}
        \includegraphics[width=\textwidth]{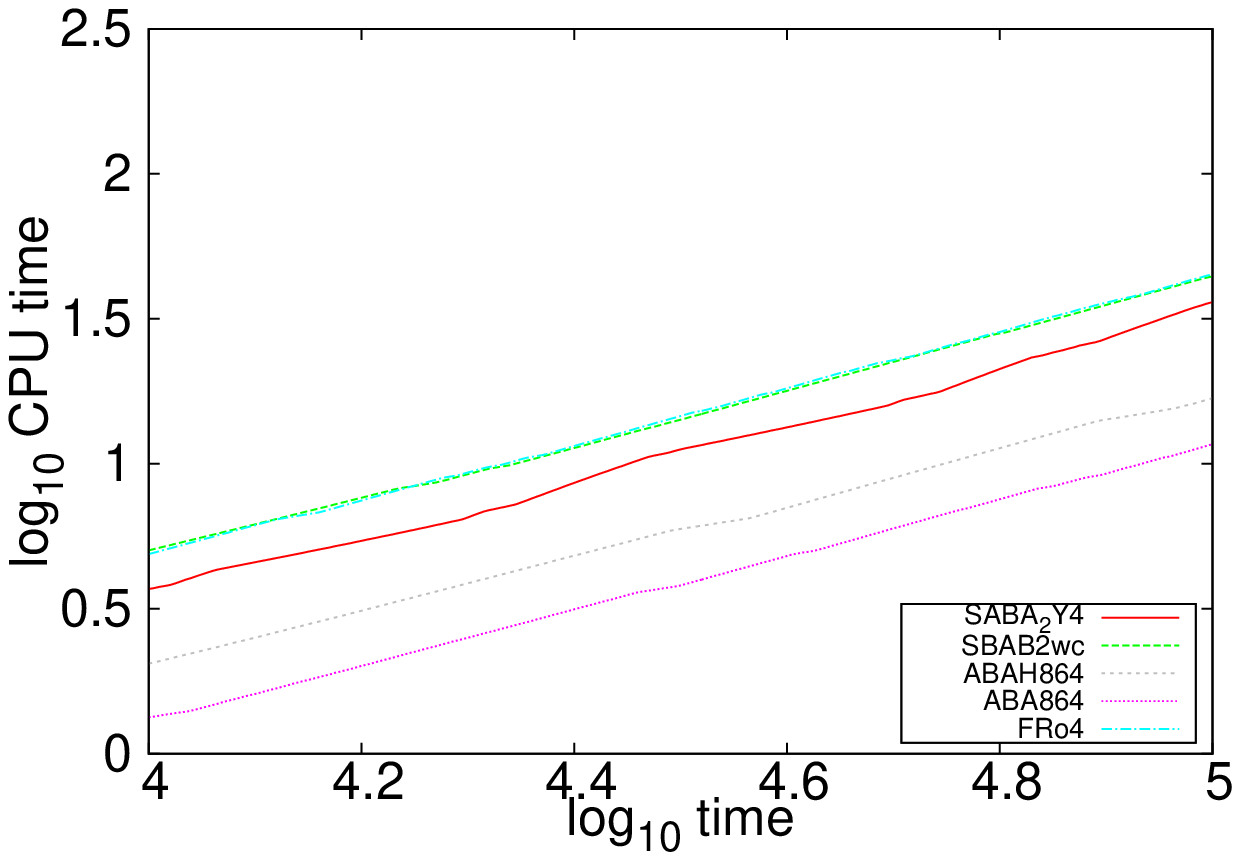}
        
    \end{subfigure}

\caption{{\tiny {\bf Results for the integration of (\ref{kgh}), by the order four schemes $SABA_2Y4$ for 
$\tau = 0.1255$ (red curves), $SBAB_2wc$ for $\tau = 0.134$ (green curve), $ABAH864$ for $\tau = 0.355$ 
(gray curve), $ABA864$ for $\tau = 0.4855$ (pink curve) and $FRo4$ for $\tau = 0.084$ 
(light blue curves).  The panels show the logarithms of the relative energy error, second moment, participation number 
and the CPU time required for evolution upto time $10^5$.}}}
\end{figure}

\newpage
 \item[{\bf 4.}] {\bf Summary and conclusions}\newline
 In this work we have studied the integration of the Klein-Gordon lattice model for the 
 so called weak chaos regime. We have used SIs of order two, four and or generalised order. The class of schemes of 
 generalised order have proven to perform better compared to the other schemes that have been tested in this study.
 Of the three schemes of generalised order, $ABA864$ performed better than $ABAH864$ and $ABA82$ in the integration of the KG model.

\subsection*{Acknowledgments}  I would like to thank Muni University for 
supporting his PhD work at the University of Cape Town through the ADBV HEST project and for facilitating him to attend 
the EAUMP conference at Makerere University. I am grateful for the input from my supervisor Dr. Ch. Skokos.

\end{enumerate}

\end{document}